\documentclass[10pt]{amsart}
\usepackage{amscd}

\title[A Result of Hermite and Equations of Degree 5 and 6]{A Result  of Hermite \\ and Equations of Degree 5 and 6}
\author{Hanspeter Kraft}
\address{Mathematisches Institut der  Universit\"at Basel,
\newline\indent
Rheinsprung 21, CH-4051 Basel, Switzerland}
\email{Hanspeter.Kraft@unibas.ch}
\date{October 2004}
\thanks{The author is supported by the Swiss National Science Foundation (Schweizerischer National\-fonds)}

\newtheorem{thm}{Theorem}
\newtheorem*{thm*}{Theorem}
\newtheorem{prop}{Proposition}
\newtheorem{lem}{Lemma}
\newtheorem*{mthm}{Main Theorem}
\newtheorem*{thA}{Theorem A}
\newtheorem*{thB}{Theorem B}

\newtheorem*{cor*}{Corollary}
\newtheorem*{conj*}{Conjecture}
\newtheorem*{prob}{Problem}

\theoremstyle{definition}
\newtheorem{defn}{Definition}
\newtheorem{exa}{Example}
\newtheorem*{exa*}{Example}

\theoremstyle{remark}
\newtheorem{rem}{Remark}

\newcommand{\op}{\operatorname}

\newcommand{\sgn}{\op{sign}}
\newcommand{\codim}{\op{codim}}
\newcommand{\qqed}{\hfill\hfill\qed}
\newcommand{\Tr}{\op{Tr}}
\newcommand{\irr}{\op{Irr}}

\newcommand{\name}[1]{\textsc{#1\/}}
\newcommand{\NN}{{\mathbb N}}
\newcommand{\ZZ}{{\mathbb Z}}
\newcommand{\QQ}{{\mathbb Q}}
\newcommand{\PP}{{\mathbb P}}
\newcommand{\FF}{{\mathbb F}}

\renewcommand{\AA}{{\mathbb A}}

\newcommand{\MMM}{\mathcal M}

\newcommand{\simto}{\xrightarrow{\sim}}
\newcommand{\be}{\begin{enumerate}}
\newcommand{\ee}{\end{enumerate}}
\newcommand{\sign}{\op{sign}}
\newcommand{\Id}{\op{Id}}
\newcommand{\gal}{\op{Gal}}
\newcommand{\SLtwo}{{\op{SL}_2}}
\newcommand{\SL}{\op{SL}}

\newcommand{\PGL}{\op{PGL}}
\newcommand{\Ind}{\op{Ind}}
\newcommand{\cha}{\op{char}}
\newcommand{\tr}{\op{tr}}

\renewcommand{\phi}{\varphi}
\begin{document}
\begin{abstract} A classical result from 1861 due to \name{Hermite} says that every separable equation of degree 5 can be transformed into an equation of the form $x^5 + b x^3 + c x + d = 0$. Later, in 1867, this was generalized to equations of degree 6 by \name{Joubert}. We show that both results can be understood as an explicit analysis of certain covariants of the symmetric groups $S_5$ and $S_6$. In case of degree 5, the classical invariant theory of binary forms of degree 5 comes into play whereas in degree 6 the existence of an outer automorphism of $S_6$ plays an essential r\^ole.
\end{abstract}

\maketitle

\section{Introduction}\label{intro}

Let $L/K$ be a finite separable field extension of degree $n$. A classical problem  is to find a generating element $x\in L$ whose equation
\begin{equation}\label{generaleq}
x^n + a_1x^{n-1} + a_2 x^{n-2} + \cdots + a_{n-1} x + a_n=0
\end{equation}
is as simple as possible. For example, a quadratic extension in characteristic $\neq 2$ has always a generator with equation $x^2 + b = 0$. One easily shows that for a separable extension of degree 3 there is a generator with equation $x^3 + b x + b = 0$ and similarly for an extension of degree 4  (see section~\ref{smalldegree}).

In degree 5 and 6 we have the following classical results which go back to \name{Hermite}  \cite{Her61} and \name{Joubert} \cite{Jou67}.
\begin{mthm}\label{mainthm} 
{\rm(a)} For any separable field extension $L/K$ of degree 5  there is a generator $x\in L$ whose equation has the form
$$
x^5 + b x^3 + c x + c = 0,
$$
except for $K=\FF_2$ where the equation has the form $x^5 + x^3 +1=0$.

{\rm(b)} Let $L/K$ be a separable extension of degree 6. If $\cha K\neq 2$  then there is a generator $x\in L$ whose equation has the form
$$
x^6 + b x^4 + c x^2 + d x + d = 0.
$$
\end{mthm}
The arguments given by \name{Hermite} and \name{Joubert} work in characteristic zero. They are  short and elegant and both are based on some classical invariant theory. The idea is to construct ``universal''  \name{Tschirnhaus} transformations which, applied to any generator of $L/K$,  produce elements of $L$ whose equations $(\ref{generaleq})$ satisfy the required properties, i.e.  $a_1=a_3 = 0$. From this it is not difficult to obtain the Main Theorem (at least in characteristic zero) although this is not explicitly formulated in their papers. To get the result also in positive characteristic and, in particular, for finite fields needs a little more work and some explicit computations.

The aim of this note is  to give a ``modern''  approach to these results, following  (and explaining)  the classical ideas.  We will show  that  \name{Hermite}'s and \name{Joubert}'s method can be considered as a careful and explicit analysis of certain covariants of the symmetric groups $S_5$ and $S_6$. In degree~5 the classical invariant theory of binary forms of degree~5 comes into play whereas in degree~6 the existence of an outer automorphism of $S_6$ plays an essential r\^ole. Another modern approach was given by  \name{Coray} in \cite{Cor87}; it is based on rationality questions for cubic hypersurfaces. 

There is the obvious question to generalize these results to higher degree. However, it was shown by \name{Reichstein} in  \cite{Rei99} that, in general, this is not possible (see Example~\ref{reichstein} at the end of section~\ref{equations}).

It should be pointed out here that the relation between equations of degree $n$ and covariants of the symmetric group $S_n$ has be studied in detail by \name{Buhler} and  \name{Reichstein} in \cite{BuR97} (see also \cite{BuR99}).

\medskip
\noindent
{\small
{\bf Acknowledgement.} I would like to thank Zinovy Reichstein and Jean-Pierre Serre for helful discussions about the subject of this paper.}

\vskip1cm
\section{Covariants and Tschirnhaus transformations}\label{tschirnhaus}

Let $K$ be field and $n\in \NN$ a positive integer. To any $x =(x_1,x_2,\ldots,x_n)\in K^n$ we associate the polynomial $\pi(x):= \prod_i (X - x_i)\in K[X]$. This defines a polynomial map
$$
\pi\colon \AA^n \to P_n
$$
where $P_n$ denotes the unitary polynomials of degree $n$: 
$$
P_n(K) := \{f = X^n + a_1 X^{n-1} + a_2 X^{n-2} + \cdots + a_{n-1} X + a_n \mid a_i \in K\}.
$$
The morphism $\pi$ is defined over  $\ZZ$ and corresponds to the {\it algebraic quotient\/} with respect to the natural action of the symmetric group $S_n$ on $\AA^n$ by permutations. This means that the polynomial functions on $P_n$ are identified, via $\pi$, with the symmetric functions on $\AA^n$,
$$
\pi^*\colon \ZZ[P_n] = \ZZ[a_1,a_2,\ldots,a_n] \simto \ZZ[\AA^n]^{S_n}, \quad a_k \mapsto 
(-1)^{k}s_k(x_1,x_2,\ldots,x_n), 
$$
where $s_k$ denotes the $k$th elementary symmetric function. If $f\in P_n(K)$ and $\xi_1,\ldots,\xi_n$ the roots of $f$ (with multiplicities) in some field extension $L/K$ then $\pi^{-1}(f)$ is the $S_n$-orbit of $\xi=(\xi_1,\ldots,\xi_n)\in \AA^n(L) = L^n$.

\medskip
Let $\Phi=(\phi_1,\phi_2,\ldots,\phi_n)\colon \AA^n \to \AA^n$ be an $S_n$-equivariant polynomial map. By definition, $\Phi$  induces a morphism $\bar\Phi\colon P_n \to P_n$ such that  the following  diagram commutes:
\[
\begin{CD}
\AA^n @>\Phi>> \AA^n \\
@VV{\pi}V  @VV{\pi}V\\
P_n @>{\bar\Phi}>> P_n
\end{CD}
\]
Such an $S_n$-equivariant morphism $\Phi$ is classically called a {\it covariant}. More generally, we have the following definition.
\begin{defn} Let  $V,W$ be finite dimensional $K$-representations of a finite group $G$. Then a $G$-equivariant $K$-morphism $\Phi\colon V \to W$ is called a {\it covariant of\/} $V$ {\it of type\/} $W$. Moreover, $\Phi$ is said to be {\it faithful\/} if $G$ acts faithfully on the image $\Phi(V_{\bar K})$ where $\bar K$ is the algebraic closure of $K$.
\end{defn}
Clearly, covariants $\Phi,\Psi\colon\AA^n \to \AA^n$  can be added and multiplied by invariants $p$ (i.e. by symmetric polynomials):

\begin{eqnarray*}
\Phi+\Psi &:=& (\phi_1+\psi_1,\phi_2+\psi_2,\ldots,\phi_n+\psi_n),\\
p\Phi &:=& (p\phi_1,p\phi_2,\ldots,p\phi_n).
\end{eqnarray*}

Thus the covariants form a module over the ring of invariants. Moreover, we can form the ``transvection'' 
$(\Phi,\Psi):=(\phi_1\psi_1,\phi_2\psi_2,\ldots,\phi_n\psi_n)$  of two covariants $\Phi$ and $\Psi$.
It is obtained by composing the product $\Phi\times\Psi$ with the bilinear covariant $\AA^n\times\AA^n \to \AA^n$, $(x_1,\ldots,x_n, y_1,\ldots,y_n)\mapsto(x_1y_1,\ldots,x_ny_n)$. 

All this also makes sense if $\Phi$ and $\Psi$ are covariants of type $\AA^n_{\text{sign}}$ where $\AA^n_{\text{sign}}$ denotes the standard representation of $S_n$ multiplied with the $\sign$ character.
E.g. if $\Psi\colon \AA^n \to \AA^n_{\sign}$ is a covariant (of type $\AA^n_{\sign}$) and $\Delta:=\prod_{i<j}(x_i-x_j)$ then $\Delta\Psi$ is a covariant of type $\AA^n$.
\begin{rem}\label{faithful} It is easy to see that if a covariant $\Phi=(\phi_1,\ldots,\phi_n)\colon\AA^n \to \AA^n$ is not faithful then $\phi_1 =\phi_2 =  \ldots = \phi_n$, and this polynomial is an $S_n$-invariant. For example, $\phi_1 + \phi_2 + \cdots +\phi_n = 0$ and if $\cha K$ does not divide $n$, then $\Phi$ is faithful if $\Phi\neq 0$.
\end{rem}
Assume now that the covariant $\Phi\colon \AA^n \to \AA^n$ is defined over $K$. Let $f\in P_n(K)$ and denote by $\xi=(\xi_1,\xi_2,\ldots,\xi_n)$ the roots of $f$ in some field extension $L/K$. Then $\Phi(\xi)=:(\bar \xi_1, \bar \xi_2,\ldots,\bar \xi_n)$ are the roots of the transformed polynomial $\bar f := \bar\Phi(f)$:
$$
f = \prod_{i=1}^n (X-\xi_i) \ \mapsto \ \bar f = \prod_{i=1}^n(X-\bar \xi_i).
$$
\begin{lem}\label{lemcov}
Let  $\Phi=(\phi_1,\phi_2,\ldots,\phi_n)\colon \AA^n \to \AA^n$ be a covariant defined over $K$. 
\be
\item $\phi_1$ is invariant under $S_{n-1} \subset S_n$ acting on the last $n-1$ variables, and $\phi_k = (1k) \phi_1$ where $(1k)$ denotes the transposition of $1$ and $k$. Conversely, every $S_{n-1}$-invariant polynomial $\phi_1$ defines a unique covariant $\Phi$ whose first component is $\phi_1$.
\item There is a uniquely defined polynomial $\phi=\phi(f,X) = \phi(a_1,\ldots,a_n,X)\in K[P_n][X]$ of degree $<n$ in $X$ such that 
$$
\phi_i(x_1,x_2,\ldots,x_n) = \phi(a_1,a_2,\ldots,a_n,x_i), \quad i=1,2,\ldots,n
$$ 
where $a_k := (-1)^k s_k(x_1,\ldots,x_n)$. Conversely, every such polynomial $\phi$ defines a covariant $\Phi\colon \AA^n \to \AA^n$.
\ee
\end{lem}
\begin{proof} Part (1) is clear, since $S_{n-1}\subset S_n$ is the stabilizer of $(1,0,\ldots,0)\in\AA^n$. The assertion (2) follows from (1) because 
$$
K[x_1,x_2,\ldots,x_n]^{S_{n-1}}= \bigoplus_{j=0}^{n-1} K[x_1,x_2,\ldots,x_n]^{S_{n}}\,x_1^{j}.
$$
This is well-known and even holds over $\ZZ$.
\end{proof}
A way to express this result is by saying that the covariants $\AA^n\to\AA^n$ form a {\it free module of rank $n$\/} over the invariants, with a basis given  the covariants 
$$
(x_1,\ldots,x_n)\mapsto (x_1^j,\ldots,x_n^j), \quad j=0,1,\ldots,n-1.
$$
For our purposes the following interpretation of Lemma~\ref{lemcov} will be important. 
If the polynomial 
$$
f = x^n + a_1 x^{n-1} + a_2x^{n-2} + \cdots+ a_n
$$ 
with coefficients $a_k\in K$ has a root $\xi$ in some field extension $L/K$ then
$\bar\xi:=\phi(f,\xi)$ belongs to $L$, and $\bar\xi$ is a root of the transformed polynomial $\bar f := \bar\Phi(f)$. Following the classics we therefore make the following definition.
\begin{defn}\label{def-tschirn}
The polynomial $\phi(f,X) \in K[P_n][X]$ from Lemma~\ref{lemcov}(2) will be called the  \name{Tschirnhaus} {\it transformation} associated to the covariant $\Phi$.
\end{defn}
Thus, by Lemma~\ref{lemcov},  the covariants $\Phi\colon\AA^n\to\AA^n$ correspond bijectively to \name{Tschirnhaus} transformations $\phi = \phi(a_1,\ldots,a_n,X)$ by
$$
\Phi(x_1,\ldots,x_n) = (\phi(a_1,\ldots,a_n,x_1),\ldots,\phi(a_1,\ldots,a_n,x_n))
$$ 
where $a_k := (-1)^k s_k(x_1,\ldots,x_n)$.
So the general problem can be formulated as follows. 
\begin{prob}
Given a field extension $L/K$ of degree $n$ and a generator $\xi \in L$ with equation $f(\xi) = 0$, find a covariant $\Phi\colon \AA^n \to \AA^n$ defined over $K$ such that $\bar f := \bar\Phi(f)$ is as simple as possible.
\end{prob}
Of course, the  transformed equation $\bar f$ has to be irreducible so that $\bar\xi=\phi(f,\xi)$ is again a generator of the extension $L/K$. 

\medskip
Our main goal is to prove the following two theorems which are  due to \name{Hermite} and \name{Joubert}. They imply the Main Theorem from the introduction about the special form of equations of degree 5 and 6. This will be shown in the following section~\ref{equations} for infinite fields $K$ and in section~\ref{finitefields} for finite fields $K$. The proofs of the two theorems below will be given in section~\ref{proofA} and \ref{proofB}.
\begin{thA}[\name{Hermite}] There is a homogeneous $S_5$-covariant $\Phi\colon \AA^5 \to \AA^5$ of degree 25, defined over $\ZZ$ and faithful for every field $K$, with the  property that
$$
s_1(\phi_1,\phi_2,\ldots,\phi_5) = s_3(\phi_1,\phi_2,\ldots,\phi_5) = 0.
$$
The covariant $\Phi$ has the form $\Phi = (\Psi,\Delta\Omega)$ where $\Psi\colon\AA^5 \to \AA^5$ is homogeneous of degree 9, $\Omega\colon\AA^5 \to \AA^5_{\text{sign}}$ homogeneous of degree 6 and $\Delta=\prod_{i<j}(x_i-x_j)$.
\end{thA}
\begin{rem}\label{Mathematica}
We used the computer program Mathematica to show that, for every prime $p$, we have $s_4(\phi_1,\ldots,\phi_5) \not\equiv 0 \mod p$.
\end{rem}
For the case of degree $6$ we need a slight modification. It is well-known that the group $S_6$ has an outer automorphism $\tau$ which is unique up to inner automorphisms. We will denote by $\AA_\tau^6$ the standard representation of $S_6$ twisted with $\tau$, i.e. $\sigma\cdot_\tau x := \tau(\sigma) x$. Clearly, $\AA^6_\tau$ has the same invariants and the same algebraic quotient $\pi\colon\AA^6_\tau \to P_6$ as $\AA^6$. Therefore every covariant $\Phi\colon \AA^6 \to \AA^6_\tau$ gives rise to a commutative diagram
$$
\CD
\AA^6 @>\Phi>> \AA^6_\tau \\
@VV{\pi}V  @VV{\pi}V\\
P_6 @>{\bar\Phi}>> P_6
\endCD
$$
It is easy to see that the results obtained so far carry over to this case with only minor modifications in the formulation, e.g. in Lemma~\ref{lemcov}(1) the component $\phi_1$ is invariant under 
$\tau(S_{5})\subset S_6$, and $\phi_k = (1k)\cdot_\tau \phi_1 = \tau((1k)) \phi_1$.
\begin{thB}[\name{Joubert}] There is a homogeneous $S_6$-covariant $\Phi\colon \AA^6 \to \AA^6_\tau$ of degree 18, defined over $\ZZ$ and faithful for every field $K$ of characteristic $\neq 2$, with the property that
$$
s_1(\phi_1,\ldots,\phi_6) = s_3(\phi_1,\ldots,\phi_6) = 0\text{ \ and \ }
s_5(\phi_1,\ldots,\phi_6) = \pm 2^s \Delta^6.
$$
The covariant $\Phi$ has the form $\Phi =\Delta  \Psi$ 
where $\Psi\colon\AA^6\to  (\AA^6_{\tau})_{\text{\rm sign}} $ is of degree 3 and $\Delta=\prod_{i<j}(x_i-x_j)$. 
\end{thB}
\begin{rem} It is interesting to remark that  $\Psi$ is the covariant of type $(\AA^6_\tau)_{\text{sign}}$ of lowest possible degree (see section~\ref{proofB} for more details).  This observation will provide us with a conceptual proof of  Theorem~B  (in characteristic zero), without any explicit calculations.
\end{rem}

\vskip1cm
\section{Equations of degree 5 and 6}\label{equations}

We will now use the results of the previous section to deduce the Main Theorem about the special form of equations of degree 5 and 6 for infinite fields $K$. For the case of finite fields we will need the explicit description of the covariants of \name{Hermite} and \name{Joubert}; this will be done in section~\ref{finitefields}. 
\begin{thm}\label{thm} Let $K$ be an infinite field.
\be
\item If  $L/K$ is a separable field extension of degree 5, then there is a generator $x$ of  $L/K$ with equation
$$
x^5 + b x^3 + c x + d = 0.
$$
\item If $L/K$ is a separable field extension of degree 6 and $\cha K \neq 2$, then there is a generator $x$ of $L/K$ with equation
$$
x^6 + b x^4 + c x^2 + dx + e = 0.
$$
\ee
\end{thm}
\begin{exa}\label{example}
Let $k$ be a field of characteristic $\neq 3$ containing a primitive third root of unity. Define
$$
L:=k(x_1,x_2,\ldots,x_r) \supset K:=k(x_1^3,x_2^3,\ldots,x_r^3).
$$
If $x\in L$, $x\neq 0$, then $\Tr_{L/K}(x^3) \neq 0$. 
In particular, there is no generator $x$ of $L/K$ whose equation has the form
$$
x^n + a_2 x^{n-2} + a_4 x^{n-4} + a_5 x^{n-5} + \cdots + a_n
$$
where $n := [L:K] = 3^r$. (In fact, if $x = \sum_{i_1,i_2,\ldots,i_r} a_{i_1 i_2\ldots i_r}x_1^{i_1}x_2^{i_2}\cdots x_r^{i_r}$ then $\Tr_{L/K} x = 0$ if and only if $a_{i_1 i_2\ldots i_r}=0$ whenever $i_1,i_2,\ldots,i_r \in3\NN$. From this observation the claim follows immeditely.)
This example also shows that for $n = 3^r$ there is no faithful covariant $\Phi\colon\AA^n\to\AA^n$ such that $s_1(\phi_1,\ldots,\phi_n) = s_3(\phi_1,\ldots,\phi_n)  =0$.
\end{exa}
In view of  Theorem A and B from the previous section the theorem above is an immediate consequence of the next proposition. In fact, applying the \name{Tschirnhaus} transformation of \name{Hermite} resp. \name{Joubert} to a ``general''  element of $L/K$ we get a new generator of $L/K$ whose  equation has  the form
$$
x^5 + b x^3 + c x + d = 0\qquad\text{resp.}\qquad
x^6 + b x^4 + c x^2 + dx + e = 0.
$$
In order to reduce further to the form of the equations claimed in the Main Theorem of the introduction (i.e. $d=c$ resp. $e=d$) we have to show that $c\neq 0$ resp. $d\neq 0$. This will be done in section~\ref{finitefields}.
\begin{prop}\label{propgeneric} 
Let $K$ be an infinite field and $L/K$ a separable field extension of degree $n$. If  $\Phi\colon \AA^n \to \AA^n$ is a faithful covariant defined over $K$ and $\phi$ the corresponding  \name{Tschirnhaus} transformation, then there is a generator $\xi$ of  $L/K$ such that $\bar\xi := \phi(\xi)$ also generates $L$ over $K$, i.e. if $f(x)=0$ is the equation of $\xi$ then the transformed equation $\bar f :=\bar\Phi(f)$ is again irreducible. Moreover, if $s\in K[x_1,\ldots,x_n]$ is a non-zero polynomial function on $\AA^n$ then $\xi\in L$ can be chosen in such a way that $s(\xi_1,\xi_2,\ldots,\xi_n)\neq 0$ where $(\xi_1,\ldots,\xi_n)$ are the conjugates of $\xi$.
\end{prop}
For the proof  we need the next two lemmas. (In the first one, the field $K$ is arbitrary.)
\begin{lem}\label{lem1} Let $\Phi\colon \AA^n\to \AA^n$ be a covariant defined over an arbitrary field $K$ and $f\in P_n(K)$ an irreducible separable polynomial with splitting field $L/K$ and Galois group $G=\gal(L/K)$. Then the image $\bar f = \bar\Phi(f)$ has the form $\bar f = h^m$ with an irreducible polynomial 
$h\in P_k(K)$ of degree $k={\frac{n}{m}}$  whose splitting field $\bar L\subset L$ is $G$-invariant. In particular,  $\gal(\bar L/K) \simeq G/N$ where $N := \gal(L/\bar L)$.
\end{lem}
\begin{proof} This is an exercise in Galois theory. If $f = \prod_{i=1}^n(X-x_i)$ and 
$\Phi(x_1,\ldots,x_n) = (\bar x_1,\ldots,\bar x_n)$, i.e. $\bar x_i = \phi(f,x_i)$, then $\bar f = \prod_{i=1}^n(X-\bar x_i)$. Since $\Phi$ is defined over $K$ we see that $\Phi\colon L^n \to L^n$ is $G$-equivariant. Hence, the splitting field $\bar L = K[\bar x_1,\ldots,\bar x_n]$ of $\bar f$ is $G$-invariant  and $G$ acts transitively on the set $\Lambda:=\{\bar x_1,\ldots,\bar x_n\}$. If $[K[\bar x_1]:K] = k$ then the stabilizer of $\bar x_1$ has index $k$ in $G$ and the orbit $\Lambda$ consists of $k$ elements, say $\Lambda=\{\bar x_1,\ldots\bar x_k\}$. It follows that  $h:=\prod_{i=1}^k (X-\bar x_i) \in K[X]$ is irreducible and $\bar f = h^m$ where $m=\frac{n}{k}$.
\end{proof}

\begin{exa} If $f\in P_n(K)$ has Galois group $S_n$ then either $\bar f$ is irreducible with Galois group $S_n$ or $\bar f = (X-a)^n$ with $a\in K$. In fact, if $x$ is a root of $f$ and $L$ the splitting field then $\gal(L/K[x]) \simeq S_{n-1}$ which is a maximal subgroup of $S_n$. Thus every element $y\in K[x]\setminus K$ generates $K[x]/K$.
\end{exa}

\begin{lem}\label{lem2} 
Let $K$ be an infinite field and $L/K$ a finite separable field extension of degree $n$. Then the subset 
$$
\irr_{L/K}:= \{f \in P_n(K) \mid f\text{ irreducible and } f( x)=0 \text{ for some }  x\in L\}
$$
is Zariski-dense in $P_n(\bar K)$  where $\bar K$ denotes the algebraic closure of $K$.
\end{lem}

\begin{proof}
Let $L=K( x)$ with equation $f( x) = 0$. A linear combination $y=\sum_{i=0}^{n-1} a_i x^i$ $(a_i\in K)$  is a generator for $L/K$ if and only if the powers $(1,y,y^2,\ldots,y^{n-1})$ are linearly independent over $K$. It follows that the corresponding subset 
$$
A:=\{(a_0,a_1,\ldots,a_{n-1})\in K^n\mid \sum_{i=0}^{n-1} a_i x^i \text{ generates } L/K\}
$$
is Zariski-dense in $\bar K^n$. As a consequence, if $ x_1:= x, x_2, \ldots,  x_n$ are the roots of $f$ in some splitting field $L'\supset L$ of $f$, then 
$$
B:=
\left\{\left( \sum_{i=0}^{n-1} a_i x_1^i, \sum_{i=0}^{n-1} a_i x_2^i,\ldots, \sum_{i=0}^{n-1} a_i x_n^i\right)
\mid (a_0,a_1,\ldots)\in A\right\} \subset {L'}^n
$$
is Zariski-dense in $\bar K^n$. Hence its image $\pi(B) \subset P_n(\bar K)$ is Zariski-dense, too. By construction, $\pi(B)$ is the set considered in the lemma.
\end{proof}
Now we can prove Proposition~\ref{propgeneric}.
\begin{proof}[Proof of Proposition~\ref{propgeneric}]  By Lemma~\ref{lem2} the set $\pi^{-1}(\irr_{L/K})$ is Zariski-dense in ${\bar K}^n$. Therefore, the subset
$$
I:= \{\xi = (\xi_1,\ldots,\xi_n)\in \pi^{-1}(\irr_{L/K})\mid \bar\xi:=\Phi(\xi)\text{ has trivial stabilizer in }S_n\}
$$
is Zariski-dense, too, because $\Phi$ is faithful. This means that the $\bar\xi_i$'s are all different and so, by Lemma~\ref{lem1}, the image $\bar f:=\bar\Phi(f)=\prod_i(X-\bar\xi_i)$ is irreducible. Also, since $I$ is Zariski-dense, the function $s$ does not vanish on $I$.
\end{proof}

\begin{exa}[\name{Reichstein} \cite{Rei99}]\label{reichstein}
Let $k$ be a field of characteristic $\neq 3$, containing the 3rd roots of unity. Define $L:=k[z_1,z_2,\ldots,z_r]\supset 
K:=k[z_1^3,z_2^3,\ldots,z_r^3]$. Thus $L/K$ is a Galois extension of degree $3^r$. It is easy to see that for every non-zero element $\xi\in L$ we have $\Tr_{L/K} \xi^3 \neq 0$. Thus there is no generator $\xi$ of $L/K$ whose equation satisfies $a_1=a_3=0$. 
\end{exa}

\vskip1cm
\section{Proof of Theorem A}\label{proofA}

For any covariant $\Phi=(\phi_1,\ldots,\phi_n)\colon\AA^n \to \AA^n$ the functions $s_k(\phi_1,\ldots,\phi_n)$ are symmetric, hence can be regarded as functions on $P_n(K)$, as we have already seen above. Denote by $W_n$ the vector space of binary forms of degree $n$:
$$
W_n(K) := K[X,Y]_n = \{f = \sum_{i=0}^n a_i X^{n-i}Y^i \mid a_i \in K\}.
$$
We will identify $P_n(K)$ with the binary forms $f$ with leading coefficient $a_0 =1$ by setting $Y=1$. Then every polynomial $q = q(a_1,\ldots,a_n)$ on $P_n(K)$ of degree $d$ defines, by homogenizing, a homogeneous polynomial $\tilde q(a_0,a_1,\ldots,a_n):=a_0^d\, q(\frac{a_1}{a_0}, \frac{a_2}{a_0},\ldots,\frac{a_n}{a_0})$ on $W_n(K)$ of the same degree. Thus,  from every covariant $\Phi\colon\AA^n \to \AA^n$ we obtain $n$ homogeneous functions $\tilde s_k\in K[W_n]$ defined by
$$
\tilde s_k(a_0,a_1,\ldots,a_n) := a_0^{d_k} s_k(\phi_1,\phi_2,\ldots,\phi_n), \quad k = 1,\ldots n,
$$
where $d_k$ is the degree of $s_k(\phi_1,\ldots,\phi_n)$ considered as a function on $P_n$ which is the same as the degree of $s_k$ in each variable $x_j$. Now recall that the group $\SLtwo(K)$ acts linearly on $W_n(K)$ by considering a binary form $f$ as a function on the standard representation $K^2$ of $\SLtwo(K)$: $gf(v) := f(g^{-1}v)$ for $g\in \SL_2(K)$ and $v\in K^2$. 

The basic idea of \name{Hermite} is to arrange the covariant $\Phi$ in such a way that the homogeneous polynomials $\tilde s_k$ become {\it $\SLtwo$-invariant\/} functions on $W_n$ and then to use our knowledge about $\SLtwo$-invariants and, in particular, the fact that there are no $\SLtwo$-invariants in certain degrees. 

In order to achieve this we will use the following classical result (see \cite[II.4 Satz 2.10]{Sch68}).
\begin{prop}\label{propschur} 
Assume that $K$ is algebraically closed of characteristic 0. Let $q \in K[x_1,\ldots,x_n]$ be a symmetric polynomial which is of degree $d$ in each variable and let $\tilde q = \tilde q(a_0.\ldots,a_n)  \in K[W_n]$  be the corresponding homogeneous polynomial of degree $d$.  Then $\tilde q$ is an $\SL_2$-invariant if and only if the following two conditions hold:
\be
\item[(T)] $q(x_1+t,x_2+t,\ldots,x_n+t) = q(x_1,x_2,\ldots,x_n)$ for  $t\in K$, i.e., $q$ only depends on the differences $x_i-x_j$;
\item[(R)] $n\cdot d$ is even and $(x_1x_2\cdots x_n)^d\, q(\frac{1}{x_1},\frac{1}{x_2},\ldots,\frac{1}{x_n}) =  (-1)^{\frac{nd}{2}} q(x_1,x_2,\ldots,x_n)$.
\ee
It then follows that $q$ is homogeneous of degree $\frac{nd}{2}$.
\end{prop}
\begin{proof}[Outline of proof] The group $\SLtwo(K)$ is generated by the matrices $\bmatrix 1 & t \\ &1\endbmatrix$ ($t\in K$) and $\bmatrix &i \\ i& \endbmatrix$ ($i:=\sqrt{-1}$). Therefore, a (homogeneous) function $\tilde q \in K[W_n]$ is $\SLtwo(K)$-invariant if and only if  $\tilde q(f)$ does not change under the following substitutions: 
$$
f(X,Y) \mapsto f(X-t,Y) \text{ and } f(X,Y)\mapsto f(iY,iX).
$$
Write $f(X,Y)= a_0\prod_{i=1}^n (X- x_iY)$ so that $\tilde q(f) = q(x_1,x_2,\ldots,x_n)$. Since
$$
f(X-t,Y) = a_0\prod_{i=1}^n((X-t)- x_iY) = a_0\prod_{i=1}^n(X-( x_i+t)Y)
$$
the invariance under $f(X,Y)\mapsto f(X-t,Y)$ is equivalent to (T). Since
$$
f(iY,iX) = a_0\prod_{k=1}^n(iY- x_kiX)= a_0(-i)^n x_1 x_2\cdots x_n 
\prod_{k=1}^n(X-\frac{1}{ x_k}Y)
$$
the invariance under $f(X,Y)\mapsto f(iY,iX)$ is equivalent to the condition
$$
q( x_1, x_2,\ldots, x_n) = 
(-i)^{nd}( x_1 x_2\cdots x_n)^d \cdot
q(\frac{1}{ x_1},\frac{1}{ x_2},\ldots,\frac{1}{ x_n})
$$
which is (R). Moreover, $\tilde q$ is also invariant under $f(X,Y)\mapsto f(tX,t^{-1}Y)$ for $t\in K^*$ which implies that $q$ is homogeneous of degree $\frac{nd}{2}$.
\end{proof}

Later on we will use the following easy fact: If two arbitrary homogeneous polynomials $q_1,q_2$ satisfy one of the conditions (T) or (R) then the same holds for the product $q_1q_2$.
\begin{exa} Start with $\Delta:=\prod_{i<j} (x_i-x_j)$. It is easy to see that the symmetric polynomial $\Delta^2$ satisfies the conditions (T) and (R) with $d:=2(n-1)$. The corresponding homogeneous invariant of degree $d$ is the   {\it discriminant} $D$ of a binary form of degree $n$.  The polynomial $\Delta$ itself satisfies the conditions (T) and (R), but is skew-symmetric, i.e. $\sigma \Delta = \sgn(\sigma)\cdot \Delta$ for $\sigma \in S_n$.
\end{exa}

Now we are ready to prove Theorem A from section~\ref{equations}.  In his note \cite{Her61} \name{Hermite}
considers the following polynomial in $\ZZ[x_1,x_2,\ldots,x_5]$:
\begin{multline}
\psi_1:=[(x_1 - x_2)(x_1-x_5)(x_4-x_3)+(x_1-x_3)(x_1-x_4)(x_2-x_5)]\cdot\\
\qquad\qquad[(x_1-x_2)(x_1-x_3)(x_5-x_4)+(x_1-x_4)(x_1-x_5)(x_2-x_3)]\cdot\\
[(x_1-x_2)(x_1-x_4)(x_5- x_3)+(x_1-x_3)(x_1-x_5)(x_4-x_2)].
\end{multline}
One easily checks that $\psi_1$ is symmetric in $x_2,x_3,\ldots,x_5$, hence, by Lemma~\ref{lemcov}, defines a covariant
$$
\Psi=(\psi_1,\psi_2,\ldots,\psi_5)\colon \AA^5 \to \AA^5
$$ 
of degree 9, defined over $\ZZ$, where $\psi_k:= (1k)\psi_1$. The functions $\psi_i$ obviously satisfy the condition (T) of Proposition~\ref{propschur}, and one finds
$$
x_1^3(x_1x_2\cdots x_5)^3 \, \psi_1(\frac{1}{x_1},\frac{1}{x_2},\ldots,\frac{1}{x_5}) = -\psi_1(x_1,\ldots,x_5).
$$
Therefore, the polynomial
\begin{equation}\label{phiHermite}
\phi_1 := \psi_1 \cdot  \prod_{1<i<j} (x_i - x_j)  \cdot \Delta
\end{equation}
has the properties  (T) and (R)  with $d=\deg_{x_i}\phi_1=10$, and $\phi_1$ is symmetric in $x_2,\ldots,x_5$. By construction, the corresponding covariant 
$$
\Phi=(\phi_1,\phi_2,\ldots,\phi_5)\colon \AA^5 \to \AA^5
$$ 
is again defined over $\ZZ$ and has degree $25 = 9 + 6 + 10$. We claim that $\Phi$ satisfies the properties of Theorem A. 

In fact, it follows from  Proposition~\ref{propschur}, choosing for $K$ the algebraic closure of $\QQ$, that the homogeneous polynomials $\tilde s_k\in K[W_5]$ corresponding to the symmetric functions $s_k(\phi_1,\ldots,\phi_5)$ are  $\SLtwo$-invariants of $W_5$ of degree $10k$. 
Since a symmetric polynomial which is divisible by $\Delta$ is automatically divisible by $\Delta^2$ we see that $\tilde s_1$ is divisible by the discriminant $D$ and $\tilde s_3$ by $D^2$, and we get $\deg \tilde s_1/D = 10 - 8 = 2$ and $\deg \tilde s_3/D^2 = 3\cdot 10 - 2\cdot 8 = 14$.
On the other hand, the $\SLtwo$-invariants of $W_5$ are generated by invariants $I_4, I_8, I_{12}$ and $I_{18}$ of degree 4, 8, 12 and 18 (see \cite[II.9 Satz 2.26]{Sch68}). Hence, there are no invariants in degree 2 and 14 and so $\tilde s_1 = \tilde s_3 = 0$ which proves  Theorem A.

Finally, the covariant $\Phi$ is faithful. First of all, $\Phi\not\equiv 0$ modulo $p$ for all primes $p$. In fact, one easily sees that the leading term of the polynomial $\psi_1$ is $-x_1^6x_2^3$ and so the leading term of $\phi_1$ has coefficient $\pm1$. Now the faithfulness  follows from Remark~\ref{faithful} for $\cha K \neq 5$. If $\cha K \neq 2$ and $\Phi$ were not faithful for $K$ then $\phi_1$ is an invariant and so $\phi_1\cdot\Delta^{-1} = \psi_1 \cdot  \prod_{1<i<j} (x_i - x_j)$ a semi-invariant, hence divisible by $\Delta$. This is not possible since $\psi_1$ does not vanish for $x_1=x_2$.
\qqed

\begin{rem} By construction, the covariant $\Phi$ has the form $\Phi = \Delta  (\Psi,\Omega) = (\Psi,\Delta\Omega)$ where 
$$
\Omega = (\omega_1,\ldots,\omega_5)\colon\AA^5\to\AA^5_{\text{sign}}
$$ 
is the homogeneous covariant of degree 6 defined by $\omega_1:= \prod_{1<i<j} (x_i - x_j)$  (and $\omega_k:=-(1k)\omega_1$ for $k\geq 2$). The representation $\AA^5_{\text{sign}}$ of $S_5$  contains the subrepresentation 
$$
U:= \{x=(x_1,\ldots,x_5)\in \AA^5 \mid x_1+\cdots+ x_5 =0\},
$$ 
and the image of the covariant $\Omega\colon \AA^5 \to \AA^5_{\text{sign}}$ is contained in $U$. (The last statement is clear since $\omega_1+\cdots+\omega_5$ is skew symmetric of degree 9, hence equal to $0$, because every skew symmetric polynomial is divisible by $\Delta$.) 

It is interesting to note  that $\Omega$ is the covariant of type $U$ of smallest possible degree because the representation $U$ occurs in $K[\AA^5]$  for the first time in degree 6. (In fact, $U$ is the irreducible representation corresponding to the partition $(2,1,1,1)$, and $K[\AA^5]_6$ contains the induced representation $\Ind_{S_2}^{S_5} K$ because the stabilizer of $x_1x_2^2x_3^3 \in K[\AA^5]_6$ is $S_2$.) 
\end{rem}
\begin{rem}
Using a computer program like \name{Singular} \cite{GPS01} or \name{Mathematica}
one can check directly  that $s_1(\phi_1,\ldots,\phi_5) = s_3(\phi_1,\ldots,\phi_5) = 0$. So the ingenious part of \name{Hermite}'s short note is the discovery of the functions $\psi_i$ above. In fact, his remark is the following see \cite{Her61}. He was trying to write out the invariant of degree 18 of the the binary forms of degree 5 in terms of the roots $x_1,\ldots,x_5$. Thus, he was  looking for a polynomial expression $\psi$ in the differences $(x_i-x_j)$ which satisfies the conditions (T) and (R) of Proposition~\ref{propschur} where $d=18$. He discovered that $\psi := \psi_1\psi_2\psi_3\psi_4\psi_5$ has this property, i.e. that $\psi$ can be written as a product of 5 terms where each one is invariant with respect to one of the standard subgroups  $S_4\subset S_5$. And, of course, he immediately realized that this can be used to transform and  simplify equations of degree 5.
\end{rem}

\vskip1cm
\section{Proof of Theorem B}\label{proofB}

We will give two proofs for Theorem B. The first one is more conceptual, but only works in characteristic zero. The second follows the explicit calculations given by \name{Joubert} and is valid in all characteristics $\neq 2$. 

\begin{proof}[First Proof]
Here the base field is $\QQ$. If $\lambda=(\lambda_1,\lambda_2,\ldots)$ is a partition of 6 we will denote by $V_\lambda$ the irreducible representation of $S_6$ associated to $\lambda$ (see \cite[\S4.1]{FuH91}). So $V_{(6)}$ is the trivial representation and $V_{(1,1,\ldots,1)}$ is the sign representation. It is not hard to see that twisting  $V_{(5,1)}$ with the outer automorphism $\tau$ we obtain the representation $V_{(2,2,2)}$ which is isomorphic to $V_{(3,3)}\otimes\sign$.

Let $V$ denote the standard representation of $S_6$, i.e. $V \simeq V_{(6)} \oplus V_{(5,1)}$. Then, as we just said, $V_{\tau} \supset V_{(2,2,2)}$. One easily sees that the third symmetric power $S^3V$ contains the representation $V_{(3,3)}$. In fact, all symmetric powers $S^iV$ are permutation representation. Since  the stabilizer of $e_1e_2e_3\in S^3V$ is $S_3\times S_3$ we see that  $S^3V$ contains the induced representation $\Ind_{S_3\times S_3}^{S_6} \QQ$ which contains $V_{(3,3)}$. 

It follows that
$$
V_{(3,3)} \simeq V_{(2,2,2)}\otimes \sign \subset V_\tau\otimes\sign
$$ 
which implies that there is a non-trivial covariant $\Psi\colon V \to V_{\tau}\otimes\sign$ of degree 3. Multiplying $\Psi$ with $\Delta$ we finally get a covariant 
$$
\Phi := \Delta\Psi \colon \AA^6 \to \AA^6_\tau
$$
of degree $3 + \deg \Delta=18$. We claim that $\Phi$ satisfies the properties of Theorem B. In fact,  for every $k$ the function $s_k(\phi_1,\ldots,\phi_6)=s_k(\psi_1,\ldots,\psi_6)\Delta^k$ is symmetric and so $s_{2k+1}(\psi_1,\ldots,\psi_6)$ is skew-symmetric of degree $6k+3$, hence is divisible by $\Delta$. Since $\deg \Delta = 15$ we get  $s_1 = s_3 = 0$. To see that $\Phi$ is faithful we simply remark that otherwise $\Phi=0$ because $s_1 = \phi_1+\cdots\phi_6= 0$ (see Remark~\ref{faithful}).
\end{proof}

\begin{proof}[Second Proof]
This proof needs some preparation. We will consider the elements of the symmetric group $S_6$ as permutations of the projective line
$$
\PP\FF_5 = \FF_5 \cup\{\infty\} = \{\infty, 0,1,2,3,4\}
$$ 
so that $H:=\PGL_2(\FF_5)$ becomes a subgroup of $S_6$ isomorphic to $S_5$. This subgroup is the image of the standard $S_5 \subset S_6$ under an outer automorphism $\tau$. Let $S_6$ act on the set of subsets of $\PP\FF_5$ consisting of 2 elements, and define $N\subset S_6$ to be the normalizer of the subset
$$
\MMM:= \{\{\infty,0\},\{1,4\},\{2,3\}\}.
$$
We obtain a surjective homomorphism $\rho\colon N \to S_3$ with kernel isomorphic to $(\ZZ_2)^3$ generated by the transpositions $(\infty\, 0),(1\,4),(2\,3)$. The following result is known and easy to prove.
\begin{lem}\label{lemgroup}
Set $N_0:= N \cap H$, the normalizer of $\MMM$ in $H$, and  $\eta:=\bmatrix 1&1\\&1\endbmatrix \in H$.
\be
\item $\rho(N_0) = S_3$ and $\ker \rho|_{N_0} \simeq (\ZZ_2)^2$.
\item $N_0$ is isomorphic to $S_4$.
\item $H = N_0 \cup \eta N_0 \cup \eta^2 N_0 \cup \eta^3 N_0 \cup \eta^4 N_0$.
\ee
\end{lem}
Now we can prove Theorem B.
For the polynomial functions on $\AA^6$ we use the variables $x_\infty,x_0,x_1,x_2,x_3,x_4$. Define, as in \name{Joubert}'s paper
\begin{eqnarray*}
h &:=& (x_\infty - x_4)(x_1 - x_3)(x_2-x_0) + (x_0-x_1)(x_4-x_2)(x_3-x_\infty) \\
&=& x_\infty x_0( x_2 + x_3 - x_1 - x_4) + x_1x_4(x_\infty + x_0 - x_2 - x_3)\\
&&\quad + x_2 x_3 (x_1 + x_4 - x_\infty - x_0)
\end{eqnarray*}
It is easy to see that $h$ is semi-invariant with respect to the subgroup $N$ defined above. 
For $\sigma\in N_0 \simeq S_4$ we have $\sigma h =\sign(\rho(\sigma))\cdot h =  \sign_{N_0}(\sigma)\cdot h$. Therefore, by Lemma~\ref{lemgroup}(3),  the function 
$$
h + \eta(h) + \eta^2 (h) + \eta^3 (h) + \eta^4 (h)
$$ 
is semi-invariant with respect to $H$. We claim that the coefficients of this polynomial are all $\pm 3$. In fact,
\begin{eqnarray*}
h &=& (x_1x_2x_3 + x_2x_3x_4 + x_4x_0x_1)
-(x_1x_2x_4 + x_2x_3x_0 + x_3x_4xx_1)\\
&+& x_\infty(x_0x_2 + x_3x_0 + x_4x_1) 
-x_\infty(x_0x_1 + x_2x_3 + x_4x_0)
\end{eqnarray*}
and each bracket expression is a sum of three monomials from a single orbit under the group generated by the cyclic permutation $\eta = (12340)\in H$. Denoting by $o_{ijk}$ the sum of the monomials in the orbit of $x_ix_jx_k$ under the group $\langle \eta \rangle\subset H$, e.g.
$$
o_{123} = x_1x_2x_3 + x_2x_3x_4 + x_3x_4x_0+ x_4x_0x_1+x_0x_1x_2,
$$
we see that
$$
h + \eta(h) + \eta^2 (h) + \eta^3 (h) + \eta^4 (h) = 
3 (o_{123} - o_{124} + x_\infty o_{02} -  x_\infty o_{01}).
$$ 
Thus
$$
\psi_1:=\frac{1}{3}(h + \eta(h) + \eta^2 (h) + \eta^3 (h) + \eta^4 (h))
$$
is the sum of all squarefree monomials $x_ix_jx_k$ ($i\neq j\neq k\neq i$) with coefficients $\pm 1$. Since $\psi_1$ is semi-invariant with respect to $H$, we see that 
$\phi_1 := \Delta \cdot \psi_1$
is invariant with respect to $H$ and defines, by Lemma~\ref{lemcov}, a homogeneous covariant
$$
\Phi = (\phi_1,\phi_2,\ldots,\phi_6) \colon \AA^6 \to \AA^6_\tau
$$
of degree $3 + \deg \Delta=18$, defined over $\ZZ$. The same degree argument as in the first proof shows that $s_1(\phi_1,\ldots,\phi_6) = s_3(\phi_1,\ldots,\phi_6) = 0$. Moreover,  in characteristic $\neq 2$, the polynomial $\psi_1$ is not a semi-invariant with respect to the whole group $S_6$, hence $\phi_1$ is not an $S_6$-invariant,  and so $\Phi$ is faithful  (see Remark~\ref{faithful}). 
\end{proof}
\begin{rem} We have seen above that $\psi_1$  is the sum of all squarefree monomials $x_ix_jx_k$ ($i\neq j\neq k\neq i$) with coefficients $\pm 1$. Thus, for any field $K$ of characteristic 2, we have  $\psi_1 = s_3$. Hence neither $\Psi$ nor $\Phi$ is  faithful in characteristic 2. We do not know if the Main Theorem for extensions of degree 6 also holds in characteristic 2.
\end{rem}

\vskip1cm
\section{The case of finite fields}\label{finitefields}

In this section we will show that the methods of \name{Hermite} and \name{Joubert} also work for finite fields thus completing the proof of the Main Theorem. 
For extensions of degree 6 in characteristic $\neq 2$  this will follows from what we have done in the previous section~\ref{equations} and~\ref{proofB}. Recall that \name{Joubert}'s covariant 
$$
\Phi\colon \AA^6 \to \AA^6_{\tau}
$$
has the form $\Phi = \Psi\cdot\Delta$ where $\Psi\colon \AA^6 \to (\AA^6_\tau)_{\sign}$ is of degree 3. If follows that $s_5(\psi_1,\ldots,\psi_6)$ is a semi-invariant of degree $15$, hence an integral multiple of $\Delta$. We claim that
$$
s_5(\psi_1,\ldots,\psi_6) = \pm 2^s\cdot \Delta \text{ \ \ for some }s\in\NN.
$$
In fact, if  $s_5(\psi_1,\ldots,\psi_6)\equiv 0\mod p$ for a prime $p\neq 2$, then it follows from Theorem~1 of section~\ref{equations} that  for any infinite field $K$ of characteristic $p$ and any extension $L/K$ of degree 6 there is a generator $\xi$ whose equation has the form $x^6 + a_2 x^4 + a_4 x^2 + a_6=0$. But this implies that $L$ contains a subfield $L' := K(\xi^2)$ of degree 3 over $K$ which clearly does not hold for generic extensions of degree 6.
\begin{rem} An explicit calculation shows that 
$$
s_5(\psi_1,\ldots,\psi_6) = -2^5 \cdot \Delta.
$$
\end{rem}
The next proposition shows that the covariant of \name{Joubert}, applied to any separabel irreducible polynomial of degree 6 over any field of characteristic $\neq 2$, always gives an irreducible polynomial. In particular, this proves the Main Theorem for extensions of degree 6.
\begin{prop}\label{propdegree6}
Let $K$ be any field of characteristic $\neq 2$ and let $f\in K[x]$ be an irreducible separable polynomial of degree 6.  If  $\Phi\colon \AA^6 \to \AA^6_\tau$ is the covariant constructed by \name{Joubert} then $\bar f = \bar\Phi(f)$ is irreducible. Moreover, the linear term of $\bar f$  has a non-zero coefficient.
\end{prop}
\begin{proof} By Theorem B we have $\bar f := \bar\Phi(f) = x^6 + b x^4 + c x^2 + d x + e$
If  $\xi =(\xi_1,\xi_2,\ldots,\xi_6)$ are the (distinct) roots of $f$, then $\Delta(\xi)\neq 0$ and so 
$$
d = s_5(\phi_1(\xi),\ldots,\phi_6(\xi))=\pm 2^s \Delta^6 \neq 0. 
$$
On the other hand, if $\bar f$ were reducible then, 
by Lemma~\ref{lem1},  $\bar f = h^k$ where $h$ is irreducible and $k = 2,3$ or $6$. The case $\bar f = h^2$ cannot occur since then $h$ should have the form $x^3 + a x$. In the other two cases a short calculation shows that $h$ is either $(x^2-a)$ or $x$. But then the coefficient $d$ of the linear term of $\bar f$  is zero.
\end{proof}
A similar result  does not hold for the covariant $\Phi$ of \name{Hermite}. In fact, if we start with an irreducible polynomial of the form $f(x) = x^5-a$ then $\bar\Phi(f) = 0$. (This can verified by using the explicit form (\ref{phiHermite}) of $\Phi$ given in section~\ref{proofA}; 
in fact,  $\psi_1(1,\zeta,\zeta^2,\zeta^3,\zeta^4) = 0$ for any fifth root of unity $\zeta$.)
However, we have the following result.
\begin{prop}\label{props4} 
Let $\Phi=(\phi_1,\phi_2,\ldots,\phi_5)\colon \AA^5 \to \AA^5$ be the covariant of \name{Hermite}.
\be
\item For every prime $p$ the symmetric polynomial $s_4(\phi_1,\ldots,\phi_5)\in\ZZ[x_1,\ldots,x_5]$ is non-zero  modulo $p$.
\item $s_4(\phi_1,\ldots,\phi_5)$ is divisible  by $\Delta^6$ in $\ZZ[x_1,\ldots,x_5]$,  and the quotient $S_4:= s_4(\phi_1,\ldots,\phi_5)/\Delta^6$ is homogeneous of degree 40.
\item If $L/K$ is a separabel extension of degree 5 and $\xi\in L$ a generator with equation $f(x)=0$ such that $S_4(\xi_1,\ldots,\xi_5) \neq 0$, then $\bar\Phi(f)$ is irreducible. ($\xi_1,\ldots,\xi_5$ are the conjugates of $\xi$.)
\ee
\end{prop}
\begin{proof} 
Recall the definition of \name{Hermite}'s covariant $\Phi\colon\AA^5 \to \AA^5$: 
$$
\phi_1 = \psi_1 \cdot \prod_{1<i<j}(x_i-x_j) \cdot \Delta
$$
where
\begin{multline*}
\psi_1:= [(x_1 - x_2)(x_1-x_5)(x_4-x_3)+(x_1-x_3)(x_1-x_4)(x_2-x_5)]\cdot\\
\qquad [(x_1-x_2)(x_1-x_3)(x_5-x_4)+(x_1-x_4)(x_1-x_5)(x_2-x_3)]\cdot\\
[(x_1-x_2)(x_1-x_4)(x_5- x_3)+(x_1-x_3)(x_1-x_5)(x_4-x_2)]
\end{multline*}
and  $\Delta = \prod_{i<j}(x_i-x_j)$.  Put $\tilde\psi_1:=\psi_1\cdot\prod_{1<i<j}(x_i-x_j)$. It is easy to see that the leading term of $\tilde\psi_1$ and $\tilde\psi_2$ is $\pm x_1^6x_2^5x_3^4$ and that the leading term of $\tilde\psi_4$ and $\tilde\psi_5$ is $\pm x_1^6x_2^6x_3^2x_4$. Moreover, one finds that $\psi_3(t^4,t^3,t^2,t,1)=0$. This implies that
\begin{multline*}
s_4(\tilde\psi_1,\ldots,\tilde\psi_5)(t^4,t^3,t^2,t,1) =\\
 \tilde\psi_1(t^4,\ldots,t,1)\cdot \tilde\psi_2(t^4,\ldots,t,1)\cdot
 \tilde\psi_4(t^4,\ldots ,t,1)\cdot \tilde\psi_5(t^4,\ldots,t,1)
\end{multline*}
and the leading term of this product is $\pm t^{188}$. Thus $s_4(\tilde\psi_1,\ldots,\tilde\psi_5)$ and hence $s_4(\phi_1,\ldots,\phi_5)=s_4(\tilde\psi_1,\ldots,\tilde\psi_5)\Delta^4$ is non-zero modulo $p$ for every prime $p$, proving (1).

We have $s_4(\phi_1,\ldots,\phi_5) = s_4(\tilde\psi_1,\ldots,\tilde\psi_5)\Delta^4$. In addition, $\tilde\psi_i\tilde\psi_j\tilde\psi_k$ is divisible by $\Delta$ for $i\neq j\neq k\neq i$, and so $s_4(\tilde\psi_1,\ldots,\tilde\psi_5$ is divisible by $\Delta^2$ since it is symmetric. Now (2) follows because $\deg \tilde\psi_i = 15$ and $\deg \Delta = 10$.

Finally, let $f\in K[x]$ be an irreducible separable polynomial of degree 5 with roots $\xi_1,\ldots,\xi_5\in \bar K$. If $\bar f := \bar\Phi(f)$ is reducible then $\bar f = (x-a)^5$ by Lemma~\ref{lem1}. If $\cha K\neq 5$ then $a=0$,  because $5a = s_1(\phi_1(\xi),\ldots,\phi_5(\xi))=0$. For $\cha K = 5$ we get $\bar f = x^5 - a^5$. In both cases we see that  $S_4(\xi_1,\ldots,\xi_5) = 0$ which contradicts the assumption. Thus we get (3).
\end{proof}

A crucial step in the proof of the Main Theorem for infinite fields $K$ was Proposition~\ref{propgeneric} which says that for a faithful covariant $\Phi\colon \AA^n \to \AA^n$ defined over $K$ and a separable extenison $L/K$ we can always find a generator $\xi\in L$ such that $\phi(\xi)$ is also a generator for $L/K$, or, equivalently, that $\bar\Phi(f)$ is irreducible where $f \in K[x]$ is the minimal polynomial of $\xi$.
However, if $K$ is finite it is not clear that such a $\xi\in L$ exists. One expects that this is the case if $K$ is large enough. In fact, we have the following more precise result. (For our proof we will only need the second part.)
\begin{prop}\label{propfinitefields}  
Let $K$ be a finite field and $L/K$ a separable extension of degree $n$. Let  $\Phi\colon \AA^n \to \AA^n$ be a faithful homogeneous covariant defined over $K$  and 
$$
\phi = p_0 + p_1 X + p_2 X^2 + \cdots + p_{n-1} X^{n-1}
$$
the corresponding  \name{Tschirnhaus} transformation (see Definition~\ref{def-tschirn}). If $\phi(\xi_1,\ldots,\xi_n)\in K$ for all generators $\xi$  of $L/K$ where $\xi_1,\ldots,\xi_n$ are the conjugates of $\xi$,  then 
$$
|K| \leq \min\{\deg p_j \mid j>0\text{ and }p_j\neq 0\} < \deg \Phi.
$$
Moreover, if $S\in K[x_1,\ldots,x_n]$ is a homogeneous symmetric polynomial such that $S(\xi_1,\ldots,\xi_n)=0$ for all generators $\xi$ of $L/K$ then
$$
|K| \leq \deg S.
$$
\end{prop}
\begin{proof} If $\phi(\xi) \in K$ for a generator $\xi$ of $L/K$ then $p_1(\xi) = p_2(\xi) = \cdots = p_{n-1}(\xi) = 0$ because $1,\xi,\xi^2,\ldots,\xi^{n-1}$ are linearly independent over $K$. Now fix a generator $\theta$ of $L/K$ and consider the following linear change of coordinates
$$
x_i = y_0 + y_1\theta_i + y_2\theta_i^2 + \cdots + y_{n-1}\theta_i^{n-1},\qquad i=1,2,\ldots,n,
$$
where $\theta_1:=\theta,\theta_2,\ldots,\theta_n$ are the conjugates of $\theta$. Each $p_j$ and also $S$ are transformed into homogeneous polynomials $\tilde p_j(y_0,\ldots,y_{n-1})$ and $\tilde S(y_0,\ldots,y_{n-1})$ of the same degree. In addition, $\tilde p_j$ and $\tilde S$ have their coefficients in $K$, because $p_j$ and $S$  do and are symmetric. 

If $\xi = a_0 + a_1\theta + \cdots + a_{n-1}\theta^{n-1}$ is a generator of $L/K$ then, by assumption, we have $\tilde p_j(a_0,\ldots,a_{n-1})=0$ for $j\geq 1$.  Thus each $\tilde p_j$ vanishes on $K^n \setminus F$ where $F$ is the finite union of all subspaces corresponding to proper subfields $L'\subset L$ containing $K$. The following Lenmma~\ref{lemfields} shows that $F$ is contained in a proper linear subspace of $K^n$ and so 
Lemma~\ref{lemzeroes} implies that $|K| \leq \deg p_j$ if $p_j\neq 0$, and also $|K|\leq \deg S$, hence the claim.
\end{proof}
\begin{lem}\label{lemfields}
Let $L/K$ be an extension of finite fields of degree $n>1$ and  $p_1,p_2,\ldots,p_k$ be the prime factors of $n$. Then the sum of the proper subfields $M\subset L$ containing $K$ has codimension
$$
\frac{n}{p_1p_2\cdots p_k}(p_1-1)(p_2-1)\cdots(p_k-1) \geq 1.
$$
\end{lem}
The following proof was communicated to me by \name{Mihaela Popoviciu} and \name{Jan Draisma}. 
\begin{proof} For every divisor $d$ of $n$ we denote by $L_d$ the (unique) subfield of $L$ with $[L:L_d]=d$. Then the the span of the proper subfields of $L$ containing $K$ is given by 
$$
L_{p_1} + L_{p_2} + \cdots +L_{p_k}.
$$
We can therefore assume that $K = L_{p_1p_2\cdots p_k}$, i.e. that $n$ is squarefree. We proceed by induction on the number $k$ of prime factors of $n$, the case $n = p_1$ being trivial. By relabeling the $p_i$'s we can assume that $p_k$ is not equal to the characteristic of $K$. Then we claim that
\begin{equation}\label{sum}
(L_{p_1}  + \cdots +L_{p_{k-1}})\cap L_{p_k} = 
L_{p_1}\cap L_{p_k}  + \cdots +L_{p_{k-1}}\cap L_{p_k}.
\end{equation}
The inclusion $\supseteq$ is clear. For the converse suppose that $\alpha\in L_{p_k}$ can be written as
$$
\alpha = \alpha_1+\alpha_2+\cdots+\alpha_{k-1} \quad \text{where } \alpha_i\in L_{p_i}.
$$
Let $F\colon L \to L$ be the \name{Frobenius} operator of the extension $L/L_{p_k}$ and put 
$$
H:=\frac{1}{p_k}(\Id + F + F^2 +\cdots + F^{p_k-1}).
$$
The linear operator $H$ is the projection onto the fixed points $L^F = L_{p_k}$ and stabilizes all $L_{p_i}$. Thus 
$ \alpha = H(\alpha) = H(\alpha_1)+ H(\alpha_2)+\cdots+H(\alpha_k)$ and 
$H(\alpha_i)\in L_{p_i}\cap L_{p_k}$
which proves our claim $(\ref{sum})$. Using this we get
\begin{eqnarray*}
\codim_K(L_{p_1}+\cdots+L_{p_k}) 
&=& \codim_K(L_{p_1}+\cdots+L_{p_{k-1}}) + \codim L_{p_k} \\
&&\ \ - \codim_K(L_{p_1} + \cdots +L_{p_{k-1}})\cap L_{p_k}\\
&=& \codim_K(L_{p_1}+\cdots+L_{p_{k-1}}) + \codim_K  L_{p_k} \\
&&\ \ - \codim_K(L_{p_1}\cap L_{p_k} + \cdots +L_{p_{k-1}}\cap L_{p_k}).
\end{eqnarray*}
Applying the induction hypothesis to the extensions $L/L_{p_1p_2\cdots p_{k-1}}$ and $L_{p_k}/K$ we find 
\begin{eqnarray*}
\codim_K(L_{p_1}+\cdots+L_{p_{k-1}}) &=& p_k\,(p_1-1)(p_2-1)\cdots(p_{k-1}-1),\\
\codim_K(L_{p_1}\cap L_{p_k} + \cdots +L_{p_{k-1}}\cap L_{p_k}) &=&
(p_1-1)\cdots(p_{k-1}-1) + \codim_KL_{p_k},
\end{eqnarray*}
hence
$$
\codim_K(L_{p_1}+\cdots+L_{p_k}) = 
(p_k-1)(p_1-1)(p_2-1)\cdots(p_{k-1}-1). 
$$
\end{proof}
\begin{lem}\label{lemzeroes}
Let $K$ be a finite field and $f\in K[y_0,y_1,\ldots,y_m]$ a non-zero homogeneous polynomial. If $f$ vanishes on $K^{m+1}\setminus W$ where $W$ is a proper linear subspace then $|K| \leq \deg f$.
\end{lem}
\begin{proof} By a linear change of coordinates we can assume that $W$ is contained in the hyperplane given by $y_0=0$. Then the polynomial $\bar f(y_1,y_2,\ldots,y_m) := f(1,y_1,\ldots,y_m)$ is non-zero,  has degree $\leq \deg f$ and vanishes on $K^m$. Now the claim follows by an easy induction on $m$, since a polynomial in one variable of degree $d$ has at most $d$ different roots.
\end{proof}
Now we are ready to give a proof of the Main Theorem for extensions of degree 5. If $K$ is infinit or $|K|\geq 40$ and $L/K$ an extension of degree 5 then there is a generator $\xi$ of $L/K$ such that $s_4(\phi_1(\xi),\ldots,\phi_5(\xi)\neq 0$, by Proposition~\ref{props4} (1) and (2) together with Proposition~\ref{propfinitefields}. It follows that the transformed equation $\bar f$ is irreducible (Proposition~\ref{props4}(3)) and has the form $x^5 + ax^3 + bx + c$ where $b\neq 0$ and the claim follows. 

\newpage
It remains to discuss the finite fields $K=\FF_q$ where $q\leq 37$ and $q\neq 2$ and to show that in each case there is an irreducible polynomial of degree 5 of the required form. It clearly suffices to consider the fields $\FF_q$ where  $q=2^2,2^3,2^5,3,5,7,11,13,17,19,23,29,31,37$. In all these cases there are the following irreducible polynomials of degree 5:
\begin{eqnarray*}
\FF_{2^2}&:& x^5+ax+a\text{ \ where $a\in\FF_{2^2}\setminus\FF_2$},\\
\FF_{2^3}&:& x^5+bx^3 + bx+b\text{ \ where $b^3+b^2+1=0$},\\
\FF_{2^5}&:& x^5+cx^3+x+1\text{ \ where $c^5+c^4+c^3+c^2+1=0$},\\
\FF_3\ &:& x^5-x-1,\\
\FF_5\ &:& x^5-x-1,\\
\FF_7\ &:& x^5-2x-2,\\
\FF_{11}&:& x^5-x-1,\\
\FF_{13}&:& x^5-x-1,\\
\FF_{17}&:& x^5+4x+4,\\
\FF_{19}&:& x^5+3x+3,\\
\FF_{23}&:& x^5+2x+2,\\
\FF_{29}&:& x^5-4x-4,\\
\FF_{31}&:& x^5+3x+3,\\
\FF_{37}&:& x^5-3x-3.
\end{eqnarray*}
This finishes the proof of the Main Theorem.

\vskip1cm
\section{Equations of degree $3$ and $4$}\label{smalldegree}
To complete the picture we want to describe the situation for equations of degree 3 and 4. First we have the following general result.
\begin{lem}\label{lemtracezero}
For $n>2$ 
there is a faithful covariant $\Phi=(\phi_1,\ldots,\phi_n)\colon \AA^n \to \AA^n_{\sign}$ of degree 
${n-1}\choose{2}$ such that $\phi_1+\cdots+\phi_n = 0$.
\end{lem}
\begin{proof} Define $\phi_1:=\prod_{1<i<j\leq n} (x_i-x_j)$ and $\phi_k:= -(1k)\phi_1$. Then $\Phi := (\phi_1,\ldots,\phi_n)$ is a faithful covariant of type $\AA^n_{\sign}$. Since $\phi_1+\cdots+\phi_n$ is skew symmetric of degree $< \deg \Delta$ the claim follows.
\end{proof}
\begin{prop}\label{propsmalldegree}
Let $L/K$ be a separable field extension of degree $n$ where $K$ is either infinite or  $\cha K$ is prime to  $n$.
\be
\item If $n>2$ there is a generator $x \in L$ with $\tr x = 0$.
\item If $[L:K]=3$ 
then there is a generator $x\in L$ which satisfies an equation of the form
$$
x^3 + a x + a = 0.
$$
\item If $[L:K]=4$ 
then there is a generator $x\in L$ which satisfies an equation of the form
$$
x^4 + a x^2 + b x + b = 0.
$$
\ee
\end{prop}
\begin{proof} (1) This is well-know if $\cha K$ is prime to $n$. If $K$ is infinite, then it follows from Lemma~\ref{lemtracezero} together with Proposition~\ref{propgeneric}.
\par\smallskip
(2) By (1) we can assume that there is a generator $x\in L$ with equation $x^3 + b x + c=0$. If $b\neq 0$ the claim follows by replacing $x$ by $\frac{b}{c}x$. If $b = 0$ (which can happen only if $\cha K \neq 3$) then $y := x + x^2$ satisfies the equation $y^3 + 3 c y + c - c^2 = 0$ which reduces to the previous case. 
\par\smallskip
(3) Again by (1) we can assume that there is a generator $x$ of $L/K$ such that $\tr x = 0$.
Thus $x$ satisfies an equation of the form $x^4 + b x^2 + c x + d = 0$. 
If $c\neq 0$ we are done as in (2). Otherwise, $\cha K \neq 2$ and the element $y:=\frac{b}{2} + x + x^2$ is again a generator of $L/K$. An easy calculation shows that the coefficient of the linear term of the equation of $y$ is equal to $a^2 - 4d$ which is non-zero because $x^4 + a x^2 + d$ is irreducible, by assumption.
\end{proof}
\begin{rem} Replacing in (2) the element $x$ by $\frac{1}{x}$ we see that for a separable extension $L/K$ of degree $3$ there is always a generator $x$ such that $x^3 + x^2 \in K$. This was mentioned to me by \name{David Masser}.
\end{rem}

\vskip1cm

\par\bigskip\bigskip
\end{document}